\documentclass[11pt]{amsart}
\usepackage{amssymb, amscd, latexsym, color, amscd, float}
\usepackage{url}
\usepackage[final]{graphicx}

\usepackage{ulem,ifthen}
\newcommand{\ifemptythenelse}[3]{%
  \begingroup
    \def\dummy{#1}%
    \def\empty{}%
    \ifx\dummy\empty{#2}\else{#3}\fi
  \endgroup
  }
\DeclareRobustCommand{\change}[2][]{%
  \ifemptythenelse{#1}{%
    \ifemptythenelse{#2}{}{\begin{color}{blue}{{#2}}\end{color}}
  }{%
    \ifemptythenelse{#2}{%
      \sout{\begin{color}{red}{{#1}}\end{color}}%
    }{%
      \ifmmode
        \begin{color}{green}{{#2}}\end{color}%
      \else
      \begin{color}{blue}{{#2}}\end{color}%
      \footnote{was: \begin{color}{red}{{#1}}\end{color}}%
      \fi
    }%
  }%
}

\newcommand{\washere}[1]{}

\newtheorem{thm}{Theorem}[section]
\newtheorem{prop}[thm]{Proposition}
\newtheorem{lem}[thm]{Lemma}
\newtheorem{cor}[thm]{Corollary}

\newtheorem{defn}[thm]{Definition}
\newtheorem{rmk}[thm]{Remark}

\begin{document}

\vspace*{2cm}


\title[\large
From local to global conjugacy]
{\large From local to global conjugacy of subgroups of relatively hyperbolic groups}

\author{Oleg Bogopolski}
\address{{Sobolev Institute of Mathematics of Siberian Branch of Russian Academy
of Sciences, Novosibirsk, Russia}\newline {and D\"{u}sseldorf University, Germany}}
\email{Oleg$\_$Bogopolski@yahoo.com}

\author{Kai-Uwe Bux}
\address{Bielefeld University, Germany}
\email{bux$\_$2009@kubux.net}


\thanks{$\dagger$ 
This research was partially supported by
SFB 701, ``Spectral Structures and Topological Methods in Mathematics'', at Bielefeld University.}


\begin{abstract}
Suppose that a finitely generated group $G$ is hyperbolic relative to a collection
of subgroups $\mathbb{P}=\{P_1,\dots,P_m\}$. Let $H_1,H_2$ be subgroups of $G$ such that $H_1$ is relatively quasiconvex with respect to $\mathbb{P}$ and $H_2$ is not parabolic. Suppose that $H_2$ is elementwise conjugate into $H_1$. Then there exists a finite index subgroup of $H_2$ which is conjugate into~$H_1$.
The minimal length of the conjugator can be estimated.

In the case where $G$ is a limit group, it is sufficient to assume only that $H_1$ is a finitely generated and $H_2$
is an arbitrary subgroup of $G$.
\end{abstract}

\maketitle

\section{Introduction}

Serre \cite[Corollary~3, page~65]{trees} famously observed that a finitely generated group acting on a tree has a global fixed point provided each of its elements fixes a point. If the tree corresponds to an {\small HNN}-extension $G=\langle B,t\rangle$, it follows that a finitely generated subgroup $H$ of $G$ conjugates into $B$ if each of its elements conjugates into $B$. If the tree corresponds to a free product $G=A * B$, the situation is similar: if each element of $H$ conjugates into $B$, then the whole subgroup conjugates into $B$. With amalgamation, however, the sutiation is not as clear cut: if each element of a finitely generated subgroup $H$ conjugates into $B$, the whole group conjugates into $B$ or into $A$. In some situations, one is able to use additional information to conclude that $H$ is conjugate into $B$. This happened to the authors in the proof of \cite[Theorem~7.4]{BogBux14}, which in turn is a central step in the argument of that paper that surface groups are subgroup conjugacy seperable ({\small SCS}). Trying to generalize, it is therefore natural to search for sufficient conditions that imply conjugacy of a subgroup into a target given that all elements conjugate into the target. Serre's observation also hints at the importance of negative curvature for this phenomenon.

In \cite{BogBux16a}, we made the first step in this direction. At the heart of the argument there lies the following:
{\medskip\par\noindent\bf Theorem~A of \cite{BogBux16a}.} {\it
  Let $G$ be a hyperbolic group, let $H_1$ be a quasiconvex subgroup of $G$ and let $H_2$ be an arbitrary subgroup of $G$. Suppose that $H_2$ is elementwise conjugate into $H_1$. Then there exists a finite index subgroup of $H_2$ which is conjugate into $H_1$.
}\medskip

We recently learned that Theorem~A was already contained in \cite[Theorem~1]{Minasyan}. We used Theorem~A to extend subgroup conjugacy separability to a class of groups including surface groups (see~\cite{BogBux16a}). Meanwhile, Chagas and Zalesskii \cite{CZ3} had generalized our result about surface groups in a different direction showing that limit groups are subgroup conjugacy separable. They use different methods.

Here, we extend Theorem~A to relatively hyperbolic groups. From a geometric point of view, a finitely generated group is hyperbolic if it acts properly discontinously and cocompactly by isometries on a hyperbolic space. In particular, point stabilizers of the action are finite. If we have an action with infinite point stabilizers, non-negative curvature might be hidden there. We would like to regard a group acting cocompactly on a hyperbolic space still as hyperbolic modulo the pieces of the group acting trivially. The notion of a group hyperbolic relative to a collection of subgroups (think of a set of representatives for point stabilizers) is one way to make this intuition precise. We recall a definition suitable for out needs and some further results in Section~\ref{relhyp}.



Our main result is:

{\medskip\par\noindent\bf Theorem~\ref{main}.} {\it
Suppose that a finitely generated group $G$ is hyperbolic relative to a collection
of subgroups $\mathbb{P}=\{P_1,\dots,P_m\}$. Let $H_1,H_2$ be subgroups of $G$ such that
$H_1$ is relatively quasiconvex with respect to $\mathbb{P}$ and $H_2$ is elementwise conjugate into~$H_1$.
Then one of the following holds:

\begin{enumerate}
\item[(1)] The subgroup $H_2$ is parabolic, i.e., conjugate into some element of $\mathbb{P}$.

\item[(2)] Some finite index subgroup of $H_2$ is conjugate into $H_1$.
\end{enumerate}

If $H_2$ is infinite and nonparabolic, then
the length of the conjugator with respect to a finite generating system $X$ of $G$ can be bounded in terms of $|X|$,
the quasiconvexity constant of $H_1$, and the minimal $X$-length of loxodromic elements of $H_2$.
}



\medskip
If all peripheral subgroups are virtually abelian, then the second option always occurs:

{\medskip\par\noindent\bf Corollary~\ref{main_cor}.} {\it
Suppose that a finitely generated group $G$ is hyperbolic relative to a collection
of virtually abelian subgroups $\mathbb{P}=\{P_1,\dots,P_m\}$. Let $H_1,H_2$ be subgroups of $G$ such that
$H_1$ is relatively quasiconvex with respect to $\mathbb{P}$ and $H_2$ is elementwise conjugate into~$H_1$.
Then some finite index subgroup of $H_2$ is conjugate into $H_1$.

Moreover, the following holds:

\begin{enumerate}

\item[(a)] If $H_2$ is infinite and nonparabolic, then
the length of the conjugator with respect to a finite generating system $X$ of $G$ can be bounded in terms of $|X|$, the quasiconvexity constant of $H_1$, and the minimal $X$-length of loxodromic elements of $H_2$.


\item[(b)] If $H_2$ is infinite and parabolic, then 
a conjugator may be chosen whose $X$-length is bounded a priory in terms of $H_1$ 
and the minimal $X$-length of elements of infinite order in $H_2$.

\item[(c)] If $H_2$ is finite then the finite index subgroup and the conjugator can be taken to be trivial.
\end{enumerate}
}

\begin{rmk}
{\rm
In Theorem~\ref{main}, passage to a finite index subgroup of $H_2$ cannot be avoided even in case that $G$ is a free group.
Indeed, let $A_4$ be the alternating group of degree 4, $K$ the Klein subgroup of $A_4$, and
$\mathbb{Z}_2$ a subgroup of $K$ of order~2.
Let $G$ be the free group of rank 2 and  $\varphi:G\rightarrow A_4$ an epimorphism.
We set $H_1=\varphi^{-1}(\mathbb{Z}_2)$ and $H_2=\varphi^{-1}(K)$. Then $H_2$ is locally, but not globally conjugate into $H_2$.}
\end{rmk}

For a subset $M$ of a relatively hyperbolic group $G$, let $M^0$ denote the subset of all loxodromic elements of $M$.

Limit groups are hyperbolic relative to a collection of representatives of
conjugacy classes of maximal noncyclic abelian subgroups~\cite[Theorem~4.5]{Dahmani}, see also~\cite[Corollary~3.5]{Alibegovich}.
This allows us to apply the main result in this case. In this situation, we also prove that the index depends only on $H_1$:

{\medskip\par\noindent\bf Corollary~\ref{MainCor}.} {\it
Let $G$ be a limit group and let
$H_1$ and $H_2$ be subgroups of $G$, where $H_1$ is finitely generated. Suppose that $H_2$ is elementwise conjugate into $H_1$. Then there exists a finite index subgroup of $H_2$ which is conjugate into~$H_1$.

The index depends only on $H_1$.
The length of the conjugator with respect to a fixed generating system $X$ of $G$
depends only on $H_1$ and $m$, where
$$m=
\begin{cases}
\underset{g\in H_2^0}\min\, dist_X(1,g) & {\rm if}\hspace*{2mm}  H_2^0\neq \emptyset.\\
\underset{g\in H_2\setminus \{1\}}\min\, dist_X(1,g) & {\rm otherwise}.
\end{cases}
$$


}

As already mentioned, our main result generalizes Theorem~A of \cite{BogBux16a}. Analogously to~\cite{BogBux16a}, our main result enables us to prove that a large class of relatively hyperbolic groups are quasiconvex-{\small SCS} and quasiconvex-{\small SICS}. In particular, we obtain a new proof that limit groups are {\small SCS}.
We will explain this in a forthcoming paper.

\section{Relatively hyperbolic groups}\label{relhyp}

The following discussion is based on \cite[Chapter~2]{Osin1}. We shall choose among the many equivalent characterizations of relatively hyperbolic groups a definition that is well suited to our needs. The first paragraphs are taken almost verbatim from \cite[Chapter~2]{Osin1}, but notation is slightly changed to suit our current needs.

Let $G$ be a group, $\{P_{\lambda}\}_{\lambda\in \Lambda}$ a collection of subgroups of $G$, $X$ a subset of $G$.
We say that $X$ is a {\it relative generating set of $G$ with respect to $\{P_{\lambda}\}_{\lambda\in \Lambda}$}
if $G$ is generated by the set $\Bigl(\underset{\lambda\in \Lambda}{\cup} P_{\lambda}\Bigr)\cup X$. (We always assume that $X$ is symmetrized, i.e. $X^{-1}=X$.) In this situation the group $G$ can be regarded as the quotient group of the free product
$$F=\bigl(\underset{\lambda\in \Lambda}\ast P_{\lambda}\bigr)\ast F(X),\eqno{(2.1)}$$
where $F(X)$ is the free group with the basis $X$.
If the kernel of the natural homomorphism $F\rightarrow G$ coincides with the normal closure of a subset $\mathcal{R}\subseteq F$,
we say that $G$ has {\it relative presentation}
$$\langle X, \{P_{\lambda}\}_{\lambda\in \Lambda}\,|\, R=1\, (R\in \mathcal{R})\rangle.\eqno{(2.2)}$$
If $X$ and $\mathcal{R}$ are finite, this relative presentation is called finite and the group $G$
is called {\it finitely presented with respect to $\{P_{\lambda}\}_{\lambda\in \Lambda}$}.

\begin{thm} {\rm (\cite[Theorem 1.1]{Osin1})}\label{finitely_presented}
Let $G$ be a finitely generated group, $\mathbb{P}=\{P\}_{\lambda\in \Lambda}$
a collection of subgroups of $G$.
Suppose that $G$ is finitely presented with respect to~$\mathbb{P}$.
Then $\Lambda$ is finite and each subgroup $P_{\lambda}$ is finitely generated.
\end{thm}

Suppose that (2.2) is a finite relative presentation. Let
$\mathcal{P}=\underset{\lambda\in \Lambda}{\sqcup}(P_{\lambda}\setminus \{1\}).$
Given a word $W$ in the alphabet $X\cup \mathcal{P}$ such that $W$ represents 1 in $G$,
there exists an expression
$$W\underset{F}{=}\,\overset{k}{\underset{i=1}{\prod}}f_i^{-1}R_if_i\eqno{(2.3)}$$
with the equality in the group $F$, where $R_i\in \mathcal{R}\cup \mathcal{R}^{-1}$ and $f_i\in F$ for any $i$.
The smallest possible number $k$ in a representation of type (2.3) is denoted $Area^{rel}(W)$.

A function $f:\mathbb{N}\rightarrow \mathbb{N}$ is called a {\it relative isoperimetric function} of (2.2)
if for any $n\in \mathbb{N}$ and for any word $W$ over $X\cup \mathcal{P}$ of length $|W|\leqslant n$
representing the trivial element of the group $G$, we have $Area^{rel}(W)\leqslant f(n)$. The smallest relative isoperimetric
function of (2.2) is called the {\it relative Dehn function} of $G$ with respect to
$\{P_{\lambda}\}_{\lambda\in \Lambda}$ and is denoted by $\delta^{rel}_{G,\{P_{\lambda},\lambda\in \Lambda\}}$
(or simply by $\delta^{rel}$ when the group $G$ and the collection of subgroups are fixed).

There are simple examples showing that $\delta^{rel}$ is not always well-defined,
i.e. it can be infinite for certain values of the argument.
However if $\delta^{rel}$ is well-defined, it is independent of the choice of the finite relative presentation up to
the following equivalence relation. Two functions $f,g:\mathbb{N}\rightarrow \mathbb{N}$ are called {\it equivalent}
if there are positive constants $A,B,C$ such that $f(n)\leqslant Ag(Bn)+Cn$
and $g(n)\leq Af(Bn) +Cn$.

\begin{defn} {\rm \cite[Definition~2.35 and Corollary~2.54]{Osin1}
{%
We call a group $G$ {\it hyperbolic relative to a collection of subgroups} $\mathbb{P}=\{P_{\lambda}\}_{\lambda\in \Lambda}$ if $G$ is finitely presented with respect to $\mathbb{P}$, the corresponding Dehn function is well-defined, and the Cayley graph $\Gamma(G,X\cup \mathcal{P})$ is a hyperbolic metric space.}
We call $\mathbb{P}$ a {\it peripheral structure} for $G$.\break
A subgroup $H\leqslant G$ is called {\it parabolic} if it is conjugate into some $P_{\lambda}\in \mathbb{P}$.

In particular, a group is hyperbolic (in the ordinary non-relative sense)
if and only if it is hyperbolic relative to the trivial subgroup. }
\end{defn}

%

Here we meet the main difficulty: The space $\Gamma(G, X\cup \mathcal{P})$
is hyperbolic, but is not locally finite
if $\mathcal{P}$ is infinite.
Suppose that $X$ generates $G$.
Then there are two distance functions on $\Gamma(G, X\cup \mathcal{P})$, namely
$dist_{X\cup \mathcal{P}}$ and $dist_X$.
For brevity, we denote $|AB|:=dist_{X\cup \mathcal{P}}(A,B)$.
Clearly, $|AB|\leqslant dist_X(A,B)$.

For reference, we collect a few statements from \cite{Arzhantseva,GM,Osin1,Osin2} that will allow us to deal with the relatively hyperbolic case.
From now on and to the end of this section, we will assume the following.

\medskip
\noindent
{\bf Assumption.} The group $G$ is generated by a finite set $X$ and is hyperbolic relative to a collection of subgroups {$\mathbb{P}=\{P_{\lambda}\}_{\lambda\in\Lambda}$}.

\begin{thm}
{\rm (\cite[Theorem 3.26]{Osin1})}\label{nu}
There exists a constant $\nu>0$ having the following property. Let $\Delta=pqr$ be a triangle whose sides $p,q,r$
are geodesics in $\Gamma(G,X\cup \mathcal{P})$. Then for any vertex $v$ on $p$, there exists a vertex $u$ on the union $q\cup r$ such that $$dist_X(u,v)< \nu.$$
\end{thm}

Recall that an element $g\in G$ is called {\it parabolic} if it is conjugate to an element of one of the subgroups $P_{\lambda}$, $\lambda\in \Lambda$. An element is called  {\it hyperbolic} if it is not parabolic. An element is called {\it loxodromic} if it is hyperbolic and has infinite order.

\begin{lem}\label{lang} {\rm (\cite[Corollary 4.20]{Osin1})}
For any loxodromic element $g\in G$, there exist $\lambda>0$, $\sigma\geqslant 0$ such that
$$dist_{X\cup \mathcal{P}}(1,g^n)\geqslant \lambda |n|-\sigma$$
for any $n\in \mathbb{Z}$.
\end{lem}



\begin{lem} {\rm (\cite[Corollary 4.21]{Osin1})}\label{equal}
Let $g$ be a loxodromic element in $G$. If $x^{-1}g^nx=g^m$ for some $x\in G$ and $n,m\in \mathbb{Z}$, then
$n=\pm m$.
\end{lem}

Recall that a subgroup of a group is called {\it elementary} if it contains a cyclic subgroup of finite index.

\begin{thm}\label{elem} {\rm (\cite[Theorem 4.3]{Osin2})} Every loxodromic element $g\in G$ is contained in a unique maximal elementary subgroup,
namely in $$E_G(g)=\{f\in G\,|\, f^{-1}g^nf=g^{\pm n}\hspace*{2mm}{\text{\rm for some}}\hspace*{2mm}{n = 1,2,3,\ldots} \}.$$
\end{thm}

\begin{lem}\label{loxodrom} {\rm (\cite[Lemma 3.5]{Arzhantseva})} Let $G$ be a group hyperbolic relative to the collection of subgroups
$\{P_{\lambda}\}_{\lambda\in \Lambda}$. For any $\lambda\in \Lambda$ and $a\in G\setminus P_{\lambda}$,
there exists a finite subset $\mathcal{F}=\mathcal{F}_{\lambda}(a)\subseteq P_{\lambda}$ such that if $x\in P_{\lambda}\setminus \mathcal{F}$,
then $ax$ is loxodromic.
\end{lem}

\begin{lem}\label{Groves_Manning} {\rm (\cite[Lemma 2.13]{GM})}
Let $H$ be a subgroup of a relatively hyperbolic group $G$. If $H$ is infinite and torsion,
then $H$ is parabolic.
\end{lem}

\begin{lem}\label{infinite_order}
Let $H$ be a subgroup of a relatively hyperbolic group $G$.
Then $H$ contains a loxodromic element if and only if
$H$ is infinite and nonparabolic.
\end{lem}

{\it Proof.} Suppose that $H$ is infinite and nonparabolic.
By Lemma~\ref{Groves_Manning}, $H$ contains an element $h$ of infinite order. Suppose that $h$ is not loxodromic.
Then $g^{-1}hg\in P_{\lambda}$ for some $g\in G$ and $\lambda\in \Lambda$.
Since $H$ is not parabolic, there exists $a\in H$ such that $g^{-1}ag\notin P_{\lambda}$.
By Lemma~\ref{loxodrom}, there exists a finite set $\mathcal{F}\subset P_{\lambda}$ such that $(g^{-1}ag)x$ is loxodromic for any $x\in P_{\lambda}\setminus\mathcal{F}$. Since $h$ has infinite order, some power $x=g^{-1}h^ng$ lies outside of~$\mathcal{F}$; hence the element $(g^{-1}ag)(g^{-1}h^ng)$ is loxodromic.
Thus, $ah^n$ is a loxodromic element in $H$. The converse direction is obvious. \hfill $\Box$

\begin{defn}
{\rm (\cite[Definitions 4.9, 4.11]{Osin1})

Let $G$ be a group generated by a finite set $X$, $\mathbb{P}=\{P_{\lambda}\}_{\lambda\in\Lambda}$ a collection of subgroups of $G$.

(a)
A subgroup $H$ of $G$ is called relatively quasiconvex with respect to $\mathbb{P}$ (or simply {\it relatively quasiconvex} when the collection $\mathbb{P}$ is fixed) if there exists $\epsilon>0$ such that the following condition holds.
Let $f_1,f_2$ be two elements of $H$ and $p$ an arbitrary geodesic path from $f_1$ to $f_2$
in $\Gamma(G,X\cup \mathcal{P})$. Then for any vertex $v\in p$, there exists a vertex $u\in H$ such that
$$dist_X(v,u)\leqslant \epsilon.$$

(b) A relatively quasiconvex subgroup $H$ of $G$ is called {\it strongly relatively quasiconvex} if
the intersection $H\cap P_{\lambda}^g$ is finite for all $g\in G$, $\lambda\in\Lambda$.

}
\end{defn}

In the case of a finitely generated relatively hyperbolic group these notions do not depend on a choice of a finite generating set:


\begin{prop}{\rm (\cite[Proposition 4.10]{Osin1})}
Let $G$ be a group hyperbolic relative to a collection of subgroups $\mathbb{P}=\{P_{\lambda}\}_{\lambda\in\Lambda}$
and let $H$ be a subgroup of $G$. Suppose that $X_1,X_2$ are two finite generating sets of $G$.
Then $H$ is (strongly) relatively quasiconvex with respect to $X_1$ if and only if it is (strongly) relatively
quasiconvex with respect to $X_2$.
\end{prop}

\begin{thm} {\rm (\cite[Theorem 4.19]{Osin1})}
Let $G$ be a finitely generated group hyperbolic relative to a collection
of subgroups $\mathbb{P}=\{P_{\lambda}\}_{\lambda\in\Lambda}$ and let $g$ be a hyperbolic element of $G$. Then the centralizer $C(g)$
of $g$ in $G$ is a strongly relatively quasiconvex subgroup in~$G$.
\end{thm}

\begin{lem}\label{Glaeser}
For every loxodromic element $b\in G$, there exists $\tau>0$
such that the following holds.
Let $m$ be a natural number and $p$ a geodesic segment in $\Gamma(G,X\cup \mathcal{P})$ connecting $1$ and $b^m$.
Then the Hausdorff distance (induced by the $dist_X$-metric) between the sets $p$ and $\mathcal{M}=\{b^i\,|\, 0\leqslant i\leqslant m\}$ is at most $\tau$.
\end{lem}

{\it Proof.} A path in $\Gamma(G,X\cup \mathcal{P})$ is called an $X$-path if its edges are labelled by elements of $X$. Let $[1,b]_X$ be a path of minimal length among all $X$-paths from 1 to $b$. We set $q_i=b^i\cdot [1,b]_X$ for $i\in \mathbb{Z}$ and $q=q_0q_1\dots q_{m-1}$.
The paths $p$ and $q$ have the same endpoints; the path $p$ is a geodesic and the path $q$ is a quasi-geodesic (with respect to $dist_{X\cup \mathcal{P}}$) with uniform constants (i.e., independent of the exponent $m$) by Lemma~\ref{lang}.

The paths $p$ and $q$ do not have backtracking~\cite[Definition~3.9]{Osin1} and
every vertex of $p$ and $q$ is a phase vertex~\cite[Definition~3.10]{Osin1}: in the case of $q$, this is because the path has only labels from $X$; and for $p$ this follows from the path being geodesic.
Thus, the technical hypotheses of \cite[Proposition~3.15]{Osin1} are satisfied.
By this proposition, there is a uniform constant $\epsilon>0$ such that
the Hausdorff distance between $p$ and $q$ (with respect to $dist_X$) is at most~$\epsilon$.
Then the Hausdorff distance betweeen $p$ and $\mathcal{M}$ is at most $\tau:=\epsilon+dist_X(1,b)$.
\hfill $\Box$

\medskip


We use the following result of B.H. Neumann.

\begin{thm}\cite[Lemma~4.1]{BN1}\label{Neumann}
Suppose that $(H_i)_{i\in I}$ is a finite family of subgroups of a group $G$ and $(x_i)_{i\in I}$ is a finite family of elements of $G$ with the property $G=\underset{i\in I}{\cup} H_ix_i$. Then there exists $i\in I$
such that $H_i$ has a finite index in $G$.
\end{thm}

\section{The main theorem}
For this section, we let $G$ be a group generated by a finite set $X$ that is hyperbolic relative to a fixed peripheral structure $\mathbb{P}=\{P_1,\dots ,P_m\}$. We put $\mathcal{P}=\overset{m}{\underset{i=1}{\sqcup}}(P_{\lambda}\setminus \{1\})$.
In the following proof we will fix a loxodromic element $b\in G$. To not interrupt the proof, we fix some notation and constants in advance.

\begin{enumerate}
\item[$\bullet$]
For an element $g\in G$, we denote $|g|:=dist_{X\cup \mathcal{P}}(1,g)$.
Clearly, $|g|\leqslant dist_X(1,g)$.

\item[$\bullet$] Let $\nu$ be the constant
defined in Theorem~\ref{nu}.

\item[$\bullet$] Let $\lambda,\sigma$ and $\tau$ be the constants defined in Lemma~\ref{lang} and Lemma~\ref{Glaeser} for the loxodromic element $b$.

\end{enumerate}

Enlarging $\nu,\sigma,\tau$ and $1/\lambda$, we may assume that all of them are integers.

\begin{thm}\label{main}
Suppose that a finitely generated group $G$ is hyperbolic relative to a collection
of subgroups $\mathbb{P}=\{P_1,\dots,P_m\}$. Let $H_1,H_2$ be subgroups of $G$ such that
$H_1$ is relatively quasiconvex with respect to $\mathbb{P}$ and $H_2$ is elementwise conjugate into~$H_1$.
Then one of the following holds:

\begin{enumerate}
\item[(1)] The subgroup $H_2$ is parabolic, i.e., conjugate into some element of $\mathbb{P}$.

\item[(2)] Some finite index subgroup of $H_2$ is conjugate into $H_1$.
\end{enumerate}

If $H_2$ is infinite and nonparabolic, then
the length of the conjugator with respect to a finite generating system $X$ of $G$ can be bounded in terms of $|X|$,
the quasiconvexity constant of $H_1$, and the minimal $X$-length of loxodromic elements of $H_2$.
\end{thm}

{\it Proof.} We may assume that $H_2$ is infinite. Suppose that $H_2$ is nonparabolic.
By Lemma~\ref{infinite_order}, $H_2$ contains a loxodromic element $b$.

Ultimately, we obtain the finite index subgroup of $H_2$ by means of Theorem~\ref{Neumann}. We shall cover $H_2$ by finitely many cosets of subgroups. One of the subgroups is the intersection $H_2\cap E_G(b)$. In order to cover the other elements, let us consider an arbitrary element $a\in H_2\setminus E_G(b)$.
Then for every $n\in \mathbb{N}$, there exists $z_n\in G$ such that $z_nab^nz_n^{-1}\in H_1$.
Let $x_n\in G$ be the shortest element with respect to the metric $dist_{X\cup \mathcal{P}}$ among those elements which satisfy the following property: There exist $k,l\geqslant 0$ such that $k+l=n$ and
$$x_n\cdot b^kab^l\cdot x_n^{-1}=z_nab^nz_n^{-1}.$$

We consider the path $\gamma_n$ starting at 1 with label $x_nb^kab^lx_n^{-1}$. The path $\gamma_n$ ends at a point $h_n\in H_1$.
We also consider the geodesic 6-gone $KABCDN$ in $\Gamma(G,X\cup \mathcal{P})$, where $K=1$, $A=x_n$, $B=x_nb^k$, $C=x_nb^ka$, $D=x_nb^kab^l$,
and $N=x_nb^kab^lx_n^{-1}$ (see Fig.~1).

\vspace*{-22mm}
\hspace*{-10mm}
\includegraphics[scale=0.7]{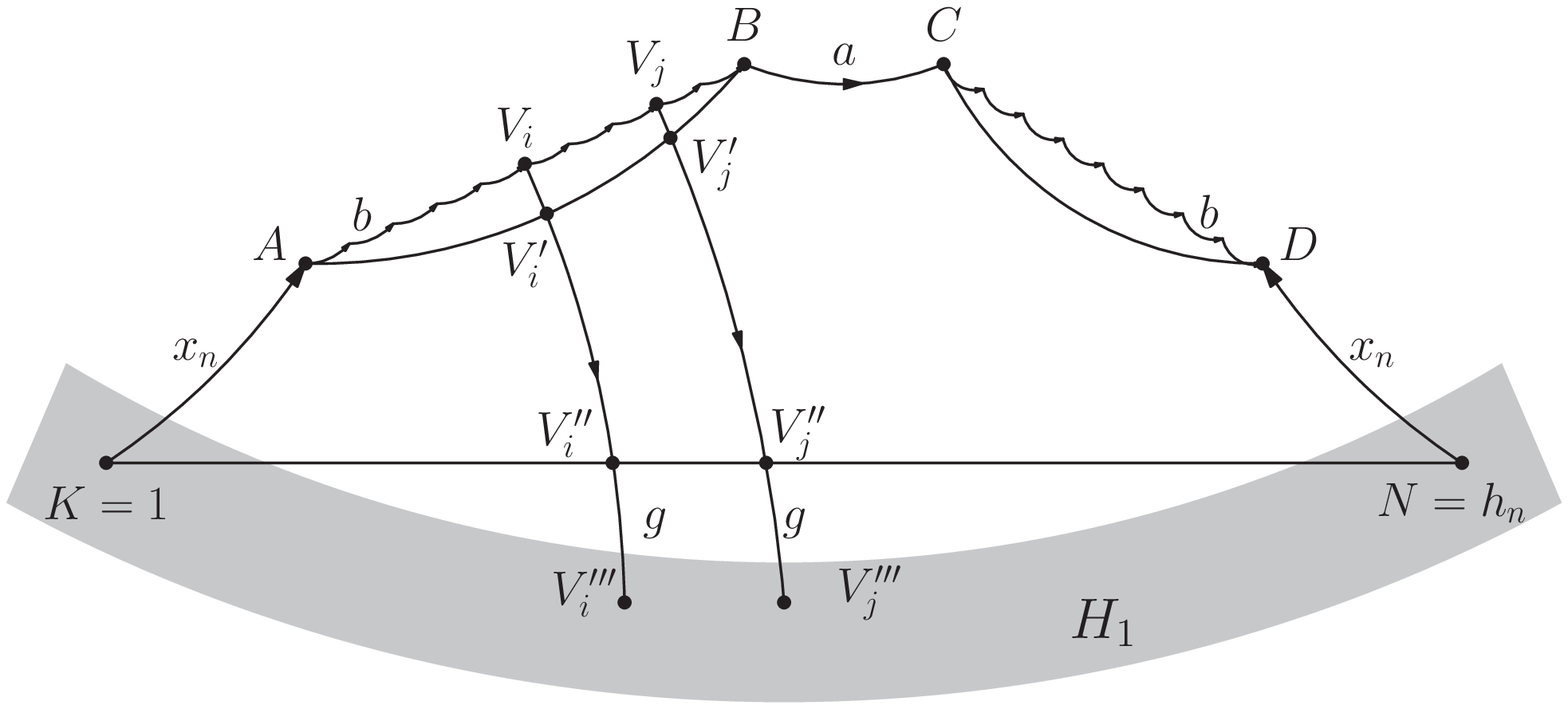}

\vspace*{-11cm}

\begin{center}
Figure 1.
\end{center}

Let $\epsilon_1$ be the quasiconvexity constant for~$H_1$. Let $\mathcal{S}(r)$ be the number of words of length up to $r$ in the alphabet $X$.
We choose
$$
n:=(\mathcal{A}+2)\cdot (\mathcal{B}_1+\mathcal{B}_2+\mathcal{B}_3+\mathcal{B}_4+2),
$$
where $\mathcal{A}: =\mathcal{S}(\tau+4\nu+\epsilon_1)$ and
$$
\begin{array}{llll}
\mathcal{B}_1: & =(\sigma+2\tau+8\nu)/\lambda,\hspace*{5mm} & \mathcal{B}_2: & =(\sigma+\tau+4\nu+|a|)/\lambda, \\
\mathcal{B}_3: & =\mathcal{S}(2\tau+4\nu),\hspace*{5mm} & \mathcal{B}_4: & =(\tau+4\nu+1)\cdot \mathcal{S}(\tau+4\nu).
\end{array}
$$

Without loss of generality, we may assume that $k\geqslant l$.
Then $k\geqslant {n/2}$.

For each $i\in [0,k]$, we consider the vertex  $V_i=A\cdot b^i$.
By Lemma~\ref{Glaeser}, there exists a vertex $V_i'\in [A,B]$ such that $$dist_X (V_i,V_i')\leqslant\tau.$$

By Theorem~\ref{nu}, there exists a point $V_i''\in [A,K]\cup [K,N]\cup [N,D]\cup [D,C]\cup [C,B]$ such that $$dist_X(V_i',V_i'')< 4\nu.$$
In particular, $$dist_X(V_i,V_i'')< 4\nu+\tau.$$

We prove that if $k$ is sufficiently large, then there are many consecutive values of $i\in [0,k]$ such that $V_i''\in [K,N]$.
We will consider four cases, where $V_i''\notin [K,N]$ and show that
each case can occur only for a restricted number of values of $i$, which we shall bound explicitly.

\medskip

{\bf Case 1.} Suppose that $V_i''\in [A,K]$.

We set $y_n:=x_nb^i$. Then $$y_n\cdot b^{k-i}ab^{l+i}\cdot y_n^{-1}=x_n\cdot b^kab^l\cdot x_n^{-1}.$$
By minimality of $|x_n|$ and using Lemma~\ref{lang}, we get:
$$
\begin{array}{ll}
|x_n|\leqslant |y_n|=|KV_i| & \leqslant |KA|-|AV_i''|+|V_i''V_i|\vspace*{2mm}\\
 & \leqslant |x_n|-|AV_i|+2|V_i''V_i|\vspace*{2mm}\\
 & \leqslant |x_n|-(i\lambda-\sigma)+2(\tau+4\nu).
\end{array}
$$
Therefore $0\leqslant i\leqslant \mathcal{B}_1$.
In particular, the number of such $i$ is at most $\mathcal{B}_1+1$.

\medskip

{\bf Case 2.} Suppose that $V_i''\in [B,C]$.

Then $$|V_iB|\leqslant |V_iV_i'|+|V_i'V_i''|+|V_i''B|\leqslant \tau+4\nu+|a|.$$
On the other side, since $V_i=x_nb^i$ and $B=x_nb^k$,  we get by Lemma~\ref{lang} that
$$|V_iB|=|b^{k-i}|\geqslant \lambda\cdot (k-i)-\sigma.$$

Thus $k-\mathcal{B}_2\leqslant i\leqslant k$.
In particular, the number of such $i$ is at most $\mathcal{B}_2+1$.

\medskip

{\bf Case 3.}
Suppose that $V_i''\in [C,D]$.

\noindent
Recall that $D=C\cdot b^l$. By Lemma~\ref{Glaeser}, there exists a vertex $W_i \in \{Cb^j\,|\, 0\leqslant j\leqslant l\}$ such that
$dist_X(V_i'',W_i)\leqslant \tau$.
Then $dist_X(V_i,W_i)\leqslant 2\tau+4\nu$.
We write $W_i=C\cdot b^{s_i}$, where $0\leqslant s_i\leqslant l$.

Since $V_i=A\cdot b^i$ and $W_i=A\cdot b^kab^{s_i}$, we have
$$dist_X(1,b^{k-i}ab^{s_i})= dist_X(V_i,W_i)\leqslant 2\tau +4\nu.$$

We claim that this case can occur for at most $\mathcal{B}_3=\mathcal{S}(2\tau+4\nu)$ values of $i\in [0,k]$.
Indeed, otherwise there would exist different $i$ and $j$ such that $b^{k-i}ab^{s_i}=b^{k-j}ab^{s_j}.$
Then $a^{-1}b^{j-i}a=b^{s_j-s_i}$. Lemma~\ref{equal} implies that $s_j-s_i=\pm (j-i)$.
By Theorem~\ref{elem}, we obtain $a\in E_G(b)$ that contradicts our assumption.

\medskip

{\bf Case 4.} Suppose that $V_i''\in [N,D]$.

We set $y_n:=x_nb^{-l}a^{-1}b^{-(k-i)}$.
Then $y_n\cdot b^{k-i}ab^{l+i}\cdot y_n^{-1}=x_n\cdot b^kab^l\cdot x_n^{-1}$.
Hence, by minimality of $|x_n|$, we have $|x_n|\leqslant |y_n|$.

Note that $y_n$ is the label of a path from $N$ to $V_i$. Therefore
$$
\begin{array}{ll}
|y_n| & =|NV_i|\leqslant |NV_i''|+|V_i''V_i|\vspace*{2mm}\\
      & =|ND|-|V_i''D|+|V_i''V_i|\vspace*{2mm}\\
      & \leqslant |x_n|-|V_i''D|+(\tau+4\nu).
\end{array}
$$

Then $|DV_i''|\leqslant \tau+4\nu$.

We claim that this case can occur for at most $\mathcal{B}_4=(\tau+4\nu+1)\mathcal{S}(\tau+4\nu)$ values of $i\in [0,k]$.
Indeed, since $V_i''$ lies on the geodesic $[D,N]$ and $|DV_i''|\leqslant \tau+4\nu$, there are at most $\tau+4\nu+1$ possibilities for $V_i''$. Since $dist_X(V_i,V_i'')\leqslant \tau+4\nu$,
there are at most $(\tau+4\nu+1)\cdot \mathcal{S}(\tau+4\nu)$ possibilities for $V_i$, and hence for $i$.


\medskip

By the choice of $n$, there exists a subset $I\subseteq [0,k]$ consisting of $\mathcal{A}+1$ consecutive integers such that
$V_i''\in [K,N]$ for every $i\in I$.
Since the subgroup $H_1$ is $\epsilon_1$-quasi\-convex, there exists
$V_i'''\in H_1$ with $dist_X(V_i'',V_i''')\leqslant \epsilon_1$. We choose an $X$-geodesic $[V_i,V_i''']$ for each $i\in I$.
We have $$dist_X(V_i,V_i''')\leqslant dist_X(V_i,V_i'')+dist_X(V_i'',V_i''')\leqslant \tau+4\nu+\epsilon_1.$$
Since $\mathcal{A}=\mathcal{S}(\tau+4\nu+\epsilon_1)$ is the number of words of length up to  $\tau+4\nu+\epsilon_1$ in the alphabet $X$,
there exist different $i,j\in I$ such that $[V_i,V_i''']$ and $[V_j,V_j''']$ have the same labels. Let $g$ be this label. Then the label of the path $[V_i'''V_i]\cup [V_i,V_j]\cup [V_j,V_j''']$ is  $g^{-1}b^{j-i}g$. We set $s_1=j-i$.
Thus, there exist $g\in G$ and $s_1\in \mathbb{Z}$ such that
$$g^{-1}b^{s_1}g\in H_1,\eqno{(3.1)}$$
$0<|s_1|\leqslant \mathcal{A}$, and $dist_X(1,g)\leqslant \tau+4\nu+\epsilon_1$.

Moreover, the label of the path $[K,A]\cup [A,V_i]\cup [V_i,V_i''']$ is $x_n b^{s_0}g$ for some $s_0$,
and we have
$$x_n b^{s_0}g\in H_1.\eqno{(3.2)}$$
Recall that $$x_n b^kab^lx_n^{-1}\in H_1.\eqno{(3.3)}$$
From (3.2) and (3.3), it follows that there are exponents $p,q\in\mathbb{Z}$ with $g^{-1}b^pab^qg\in H_1$, and using (3.1), we can in addition arrange for $0\leqslant p,q\leqslant |s_1|\leqslant~\mathcal{A}$. Thus,
$$a\in b^{-p}gH_1g^{-1}b^p\cdot b^{-(q+p)}.$$
Since $a$ is an arbitrary element in $H_2\setminus E_G(b)$ and $dist_X(1,g)\leqslant \tau+4\nu+\epsilon_1$,
we have
$$H_2\subseteq \underset{(t,z)\in M}{\bigcup}(z^{-1}H_1z\cdot b^{-t})\cup E_G(b),$$
where $M=\{(t,z)\in \mathbb{Z}\times G\,|\, 0\leqslant t\leqslant 2\mathcal{A},\, dist_X(1,z)\leqslant
\mathcal{A}\cdot dist_X(1,b)+\tau+4\nu+\epsilon_1\}.$

Since $b\in H_2$, we have
$$H_2= \underset{(t,z)\in M}{\bigcup}((z^{-1}H_1z\cap H_2)\cdot b^{-t})\cup (E_G(b)\cap H_2).$$
Since the set $M$ is finite, we deduce from Theorem~\ref{Neumann} that either $E_G(b)\cap H_2$ is of finite index in $H_2$, or there exists $(t,z)\in M$ such that
$z^{-1}H_1z\cap H_2$ is of finite index in $H_2$.
In the first case, $\langle b\rangle$ has finite index in $H_2$ and we are done.
In the second case, a finite index subgroup of $H_2$ is conjugate into $H_1$.
\hfill $\Box$




\medskip

\section{Toral relatively hyperbolic groups}

In this section we specialize Theorem~\ref{main} for the case,
where the peripheral subgroups of the relatively hyperbolic group $G$ are virtually abelian.


\begin{thm}\label{malnormality}
{\rm (\cite[Theorem 1.4]{Osin1})}
Let $G$ be a group, $\mathbb{P}=\{P_{\lambda}\}_{\lambda_\in \Lambda}$ a collection of subgroups of $G$.
Suppose that $G$ is finitely presented with respect to $\mathbb{P}$
and the relative Dehn function of $G$ with respect to $\mathbb{P}$ is well-defined, i.e., it takes finite values for each $n\in \mathbb{N}$.
Then the following conditions hold.

\begin{enumerate}
\item[{\rm (1)}] For any $g_1,g_2\in G$ the intersection $P_{\lambda}^{g_1}\cap P_{\mu}^{g_2}$ is finite
whenever $\lambda\neq \mu$.

\item[{\rm (2)}] The intersection $P_{\lambda}^g\cap P_{\lambda}$ is finite for any $g\notin P_{\lambda}$.
\end{enumerate}
\end{thm}

A quasiconvex subgroup of a relatively hyperbolic group is hyperbolic relative to an induced peripheral structure:
\begin{thm}\label{Hrushka}

{\rm (\cite[Theorem 9.1]{Hruska1})}
Let $G$ be a finitely generated group hyperbolic relative to a collection of subgroups
$\mathbb{P}=(P_{\lambda})_{\lambda\in \Lambda}$ and let $H\leqslant G$ be a relatively quasiconvex subgroup.
Consider the following collection of subgroups of $H$:
$$\overline{\mathbb{O}}:=\{H\cap gPg^{-1}\,|\, g\in G,\,\, P\in \mathbb{P},\hspace*{2mm}{\text{\it and}}\hspace*{2mm} H\cap gPg^{-1}\hspace*{2mm}{\text{\it is infinite}}\}.$$
Then the elements of $\overline{\mathbb{O}}$ lie in only finitely many conjugacy classes of $H$.
Furthemore, if $\mathbb{O}$ is a set of representatives of these classes, then $H$ is hyperbolic relative to $\mathbb{O}$. We call $\mathbb{O}$ the {\it induced peripheral structure of $H$.}
\end{thm}

\begin{lem}\label{NoLoxodromicElement}
Let $G$ be a finitely generated group hyperbolic relative to a collection
of subgroups $\mathbb{P}=(P_{\lambda})_{\lambda\in \Lambda}$.
Let $H_1,H_2$ be subgroups of $G$ such that $H_1$ is relatively quasiconvex
with respect to $\mathbb{P}$ and $H_2$ is a subgroup of a peripheral subgroup $P_{\lambda}\in \mathbb{P}$.\\
Suppose that $H_2$ is elementwise conjugate into~$H_1$.
Then there exists a finite collection of elements $g_1,\dots ,g_n\in G$
such that each element of infinite order of $H_2$ is conjugate (in $P_{\lambda}$) into the union of $P_{\lambda}\cap g_i^{-1}H_1g_i$. The elements $g_1,\dots ,g_n\in G$ can be chosen so that they depend only on $H_1$,
but not on~$H_2$.
\end{lem}

{\it Proof.} We will use notation of Theorem~\ref{Hrushka} applied to $H=H_1$.
First, we fix a peripheral structure of $H_1$ induced by the peripheral structure of $G$:
Then there exists a finite collection of elements $(g_i)_{i\in I}$ of $G$ and finite subsets $\Lambda_i\subset \Lambda$  for each $i\in I$ such that

\begin{enumerate}
\item[$\bullet$]
$\mathbb{O}=(H_1\cap g_iP_{\mu}g_i^{-1})_{i\in I, \mu\in \Lambda_i}$ consists
of representatives of conjugacy classes in $H_1$ of subgroups from $\overline{\mathbb{O}}$
and

\item[$\bullet$] $H_1$ is hyperbolic relative to $\mathbb{O}$.
\end{enumerate}

Recall that $H_2\leqslant P_{\lambda}$ for some $\lambda\in \Lambda$. Let $h\in H_2$ be an element of infinite order. We show that $h$ is conjugate (in $P_{\lambda}$) into $P_{\lambda}\cap g_i^{-1}H_1g_i$ for some $i\in I$.
By assumption, there exists $g\in G$ such that $ghg^{-1}\in H_1$.
Thus, $$ghg^{-1}\in H_1\cap gP_{\lambda}g^{-1}.\eqno{(4.1)}$$
In particular, $H_1\cap gP_{\lambda}g^{-1}$ is infinite.
By choice of $\mathbb{O}$, there exist $i\in I$, $\mu\in \Lambda_i$ and $z\in H_1$ such that
$$z^{-1}(H_1\cap gP_{\lambda}g^{-1})z=H_1\cap g_iP_{\mu}g_i^{-1}.\eqno{(4.2)}$$
It follows that $z^{-1}gP_{\lambda}g^{-1}z \cap g_iP_{\mu}g_i^{-1}$ is infinite, whence
by Theorem~\ref{malnormality}, we have $\lambda=\mu$ and $x:=g^{-1}zg_i\in P_{\lambda}$.
Now (4.1) and (4.2) imply $x^{-1}hx\in P_{\lambda}\cap g_i^{-1}H_1g_i$, and we are done. \hfill $\Box$


\medskip

\begin{lem} {\rm (\cite[Lemma~5.8]{Osin1})}\label{Kakao}
Let $G$ be a group generated by a finite set $X$.
Suppose that $G$ is hyperbolic relative to a collection
of recursively presented subgroups $\mathbb{P}=\{P_1,\dots,P_m\}$, each with solvable conjugacy problem.
There exists a computable function $f:\mathbb{N}\rightarrow \mathbb{N}$ with the following property: if $g\in G$ is conjugate into $P_{\lambda}$, then $z^{-1}gz\in P_{\lambda}$ for some $z\in G$ with $dist_X(1,z)\leqslant f(dist_X(1,g))$.
\end{lem}

Since finitely generated virtually abelian groups are finitely presented and have solvable conjugacy problem, we obtain the following:
\smallskip
\begin{cor}~\label{main_cor}
Suppose that a finitely generated group $G$ is hyperbolic relative to a collection
of virtually abelian subgroups $\mathbb{P}=\{P_1,\dots,P_m\}$. Let $H_1,H_2$ be subgroups of $G$ such that
$H_1$ is relatively quasiconvex with respect to $\mathbb{P}$ and $H_2$ is elementwise conjugate into~$H_1$.
Then some finite index subgroup of $H_2$ is conjugate into $H_1$.

Moreover, the following holds:

\begin{enumerate}

\item[(a)] If $H_2$ is infinite and nonparabolic, then
the length of the conjugator with respect to a finite generating system $X$ of $G$ can be bounded in terms of $|X|$, the quasiconvexity constant of $H_1$, and the minimal $X$-length of loxodromic elements of $H_2$.


\item[(b)] If $H_2$ is infinite and parabolic, then a conjugator may be chosen 
whose $X$-length is bounded a priory in terms of $H_1$ and the minimal $X$-length 
of elements of infinite order in $H_2$.

\item[(c)] If $H_2$ is finite then the finite index subgroup and the conjugator can be taken to be trivial.
\end{enumerate}

\end{cor}

{\it Proof.} We assume that $H_2$ is infinite, otherwise the statement is trivial.
If $H_2$ is nonparabolic, then the statement follows from Theorem~\ref{main}. Thus, we assume that $H_2$ is infinite and parabolic. Then there exist $x\in G$ and $P\in \mathbb{P}$ such that $x^{-1}H_2x\leqslant P$.

Let $g\in H_2$ be the element of infinite order with the minimal $X$-length.
Since $x^{-1}gx\in P$, we have $z^{-1}gz\in P$ for some $z\in G$ with $dist_X(1,z)\leqslant f(dist_X(1,g))$
(see Lemma~\ref{Kakao}). Then $P^{z^{-1}x}\cap P$ is infinite. By Lemma~\ref{malnormality} we have $x\in zP$. Hence, $z^{-1}H_2z\leqslant P$, where the $X$-length of $z$ is bounded as above. Thus,
we may assume that $H_2$ lies in $P$.

Let $P'$ be a maximal torsion free abelean subgroup of $P$. Clearly, $k:=|P:P'|$ is finite.
Let $z_1,\dots,z_k$ be a set of representatives of left cosets of $P'$ in $P$.
We set $H_2':=H_2\cap P'$.

By Lemma~\ref{NoLoxodromicElement},
there exists a finite collection of elements $g_1,\dots ,g_n\in G$
such that $H_2'$ is elementwise conjugate (in $P$) into the union of $P\cap H_1^{g_i}$.
Then $H_2'$ lies in the union of $P\cap H_1^{g_iz_j}$, and we have
$$H_2'=\overset{n}{\underset{i=1}{\cup}}\overset{k}{\underset{j=1}{\cup}} (H_2'\cap H_1^{g_iz_j}).$$
By Theorem~\ref{Neumann}, one of the subgroups $H_2'\cap H_1^{g_iz_j}$ has finite index in $H_2'$ and hence in $H_2$.
Recall that by Lemma~\ref{NoLoxodromicElement}, the elements $g_1,\dots,g_n$ depend only on $H_1$, and the elements $z_1,\dots,z_k$ depend only on $P$.
\hfill $\Box$

\section{An application to limit groups}

Recall that a subgroup $H$ of a group $G$ is called a {\it retract} of $G$ if there exists an epimorphism $f:G\rightarrow H$ with $f_{|H}={\rm id}$. The epimorphism $f$ is called a {\it retraction}. Equivalently, $H$ is retract of $G$ if $G=K\rtimes H$ for some subgroup $K$ of $G$. A subgroup $H$ of $G$ is called a {\it virtual retract} of $G$ if $H$ is a retract of a finite index subgroup of~$G$.
In the following proposition we list some properties of limit groups that we use later.

\begin{prop}\label{properties of limit groups}
\begin{enumerate}
\item[(1)] Limit groups have the unique root property, i.e., $x^n=y^n$ with $n\neq 0$ implies $x=y$.
In particular, limit groups are torsion free.

\item[(2)] Retracts of limit groups are closed under taking of roots.



\item[(3)] If $H$ is a finitely generated subgroup of a limit group $G$,
then $H$ is a virtual retract of $G$~\cite[Theorem B]{W}, and hence
is quasi-isometrically embedded~\cite[Corollary 3.12]{W}. In particular, $H$ is quasiconvex in $G$.
\end{enumerate}
\end{prop}
{\it Proof.}
(1) It suffices to argue the unique root property for subgroups of limit groups generated by two elements. By \cite[Proposition~3.1\,(4)]{CD}, those subgroups are free or free abelian; and those have the unique root property.

(2) Let $G$ be a limit group and $f:G\rightarrow H$ a retraction.
Suppose that $g\in G$ is such that $g^n\in H$ for some $n\in \mathbb{N}$.
Then $f(g)^n=f(g^n)=g^n$. Since $G$ has the unique root property, we have $g=f(g)\in H$.
\hfill$\Box$



\begin{cor}\label{MainCor}
Let $G$ be a limit group and let
$H_1$ and $H_2$ be subgroups of $G$, where $H_1$ is finitely generated. Suppose that $H_2$ is elementwise conjugate into $H_1$. Then there exists a finite index subgroup of $H_2$ which is conjugate into~$H_1$.

The index depends only on $H_1$.
The length of the conjugator with respect to a fixed generating system $X$ of $G$
depends only on $H_1$ and $m$, where
$$m=
\begin{cases}
\underset{g\in H_2^0}\min\, dist_X(1,g) & {\rm if}\hspace*{2mm}  H_2^0\neq \emptyset.\\
\underset{g\in H_2\setminus \{1\}}\min\, dist_X(1,g) & {\rm otherwise}.
\end{cases}
$$

\end{cor}

{\it Proof}.
Let $\mathbb{P}$ be a set of representatives for the conjugacy classes of maximal non-cyclic abelian subgroups of $G$. As a limit group, $G$ is hyperbolic relative to the finite family $\mathbb{P}$
(\cite[Theorem~4.5]{Dahmani}, another proof is given in \cite[Corollary 3.5]{Alibegovich}).

By statement (3) of Proposition~\ref{properties of limit groups},
$H_1$ is quasiconvex.
By Corollary~\ref{main_cor},
there exists a finite index subgroup $\widetilde{H_2}$ of $H_2$ and an element $g\in G$
such that $g^{-1}\widetilde{H_2}g\leqslant H_1$ and the $X$-length of $g$ is bounded as above.
It remains to prove that the index $|H_2:\widetilde{H}_2|$ only depends on $H_1$.


By statement (3) of Proposition~\ref{properties of limit groups},
there exists a finite index subgroup $\widetilde{G}$ of $G$ such that $H_1$ is a retract of $\widetilde{G}$.
Let $N$ be a finite index subgroup in $\widetilde{G}$ which is normal in $G$.

From $g^{-1}\widetilde{H_2}g\leqslant H_1$ we deduce $(g^{-1}H_2g)^{n!}\leqslant H_1$, where $n=|H_2:\widetilde{H_2}|$.
In particular, $$\bigl(g^{-1}(H_2\cap N)g\bigr)^{n!}\leqslant H_1.$$

Since $H_1$ is a retract of $\widetilde{G}$ and $g^{-1}(H_2\cap N)g\leqslant N\leqslant \widetilde{G}$,
we deduce from the statement (2) of Proposition~\ref{properties of limit groups} that $$g^{-1}(H_2\cap N)g\leqslant H_1.$$
The estimate $|H_2:H_2\cap N|\leqslant |G:N|$ completes the proof.
\hfill $\Box$

\end{document}